\documentclass[a4paper,11pt]{article}

\usepackage[top=3cm,bottom=3cm,left=2.5cm,right=2.5cm]{geometry}
\usepackage[T1]{fontenc}
\usepackage{iftex}
\ifPDFTeX
  \usepackage[utf8]{inputenc}
\fi
\usepackage{lmodern,microtype}
\usepackage[english]{babel}
\usepackage{amsmath,amssymb,amsthm,mathtools}
\usepackage{enumitem}
\usepackage{authblk}
\usepackage{cite}
\usepackage[hidelinks]{hyperref}
\hypersetup{pdftitle={Equality in the clique bound for critical graphs},
  pdfauthor={Hongyang Wang}}
\setlist[enumerate]{leftmargin=2.2em,itemsep=0.15em,topsep=0.25em}
\allowdisplaybreaks
\newtheorem{theorem}{Theorem}[section]
\newtheorem{lemma}[theorem]{Lemma}
\newtheorem{problem}[theorem]{Problem}
\newcommand{\QQ}{\mathbb{Q}}

\title{Equality in the clique bound for critical graphs}
\author{Hongyang Wang\thanks{Email address: \href{mailto:hwanghj@connect.ust.hk}{hwanghj@connect.ust.hk}}}
\affil{Department of Mathematics,\\
Hong Kong University of Science and Technology,\\
Clear Water Bay, Hong Kong}
\date{}

\begin{document}
\maketitle

\begin{abstract}
A graph is $k$-critical if its chromatic number is $k$ and every proper subgraph is $(k-1)$-colorable. Gao and Ma (2023) proved that every $n$-vertex $k$-critical graph with $n > k \ge 4$ contains at most $n - k + 3$ cliques of size $k - 1$, and asked which graphs attain this bound. We answer this question by showing that equality holds if and only if $n - k + 3$ is odd and the graph is isomorphic to the join of a complete graph on $k - 3$ vertices and a cycle of length $n - k + 3$. More generally, for $n>k\ge3$, we determine when the size of a nonempty family of distinct $(k-1)$-cliques in an $n$-vertex $k$-critical graph attains the upper bound $s-k+3$, where $s$ is the number of vertices covered by the family. Equality holds if and only if the graph is such a join (with the complete part empty when $k = 3$) and the family consists of all its $(k-1)$-cliques.
\end{abstract}



\section{Introduction}\label{sec:intro}

A finite simple graph $G$ is \emph{$k$-critical} if its chromatic number $\chi(G)$ is $k$ and every proper subgraph is $(k-1)$-colorable. The systematic study of critical graphs goes back to Dirac~\cite{Dirac1952}. An $r$-clique in $G$ is a set of $r$ pairwise adjacent vertices of $G$. Write $t_r(G)$ for the number of $r$-cliques in $G$ and $K_r$ for the complete graph on $r$ vertices. Since the $3$-critical graphs are precisely the odd cycles, every $n$-vertex $3$-critical graph $G$ satisfies $t_2(G) = n$. Gallai conjectured that $t_{k-1}(G) \le n$ holds for every $n$-vertex $k$-critical graph $G$~\cite{JensenToft1995}. For $k = 4$, this asserts that an $n$-vertex $4$-critical graph contains at most $n$ triangles; by contrast, Ma and Yang~\cite{MaYang2021} proved that every $n$-vertex $4$-critical graph contains at least a constant multiple of $n^2$ odd cycles.

Stiebitz~\cite{Stiebitz1987} proved Gallai's bound $t_{k-1}(G) \le n$ for $k = 4$, and Abbott and Zhou~\cite{AbbottZhou1992} proved it for all $k \ge 4$, showing further that equality holds only when $G \cong K_k$. Excluding $K_k$, Abbott and Zhou~\cite{AbbottZhou1992} asked whether $t_{k-1}(G) \le n - k + 3$ holds whenever $k \ge 5$ and $n \ge k + 2$. Su obtained this bound for $k = 5, 6$ in~\cite{Su1995DM} and for $k = 7$ in~\cite{Su1995JCMCC}, while K\'ezdy and Snevily~\cite{KezdySnevily1997} proved that $t_{k-1}(G) < n - 3k/5 + 2$ for all $n > k \ge 3$. Gao and Ma~\cite{GaoMa2023} settled the question in full by proving that $t_{k-1}(G) \le n - k + 3$ for every $n$-vertex $k$-critical graph $G$ with $n > k \ge 4$. Combining this clique bound with Tur\'an's theorem, Luo, Ma, and Yang~\cite{LuoMaYang2023} obtained upper bounds on the number of edges in $k$-critical graphs.

Write $C_m$ for the cycle on $m$ vertices. For vertex-disjoint graphs $H_1$ and $H_2$, let $H_1 \vee H_2$ denote their \emph{join}, obtained from the disjoint union of $H_1$ and $H_2$ by adding an edge between every vertex of $H_1$ and every vertex of $H_2$. For $k \ge 4$ and odd $m \ge 5$, the join $K_{k-3} \vee C_m$ is an $n$-vertex $k$-critical graph with $n = k + m - 3$. Because $C_m$ is triangle-free, each $(k-1)$-clique of $K_{k-3} \vee C_m$ consists of $V(K_{k-3})$ together with the endpoints of an edge of $C_m$, so $t_{k-1}(K_{k-3} \vee C_m) = m = n - k + 3$. Abbott and Zhou~\cite{AbbottZhou1992} had already shown that for $k = 4$, only the odd wheels $K_1 \vee C_m$ attain $t_3(G) = n - 1$. Gao and Ma~\cite{GaoMa2023} suggested that $K_{k-3} \vee C_{n-k+3}$ is the unique extremal graph whenever $n - k + 3$ is odd, and posed the general equality problem.

\begin{problem}[Gao and Ma~\cite{GaoMa2023}]\label{prob:gaoma}
Determine all $n$-vertex $k$-critical graphs $G$ with $n > k \ge 4$ satisfying $t_{k-1}(G) = n - k + 3$.
\end{problem}

The following theorem answers Problem~\ref{prob:gaoma}.

\begin{theorem}\label{thm:equality}
Let $G$ be an $n$-vertex $k$-critical graph with $n > k \ge 4$, and let $m = n - k + 3$. Then $t_{k-1}(G) = m$ if and only if $m$ is odd, $m \ge 5$, and $G \cong K_{k-3} \vee C_m$.
\end{theorem}

Together with the upper bound of Gao and Ma~\cite{GaoMa2023}, Theorem~\ref{thm:equality} gives $t_{k-1}(G) \le n - k + 2$ whenever $n > k \ge 4$ and $n - k + 3$ is even.

We deduce Theorem~\ref{thm:equality} from a local bound that replaces the order $n$ of $G$ with the number of vertices covered by a chosen family of cliques. For a family $\mathcal{A}$ of cliques of $G$, the \emph{support} of $\mathcal{A}$ is the vertex set $S_{\mathcal{A}} = \bigcup_{Q \in \mathcal{A}} Q$. In the next theorem we also allow $k = 3$, where $K_0$ has no vertices and $K_0 \vee C_m = C_m$.

\begin{theorem}\label{thm:selected}
Let $G$ be an $n$-vertex $k$-critical graph with $n > k \ge 3$, and let $\mathcal{A}$ be a nonempty family of distinct $(k-1)$-cliques of $G$. Then $|\mathcal{A}| \le |S_{\mathcal{A}}| - k + 3$. Equality holds if and only if $G \cong K_{k-3} \vee C_m$ for an odd integer $m = n - k + 3 \ge 5$ and $\mathcal{A}$ consists of all $(k-1)$-cliques of $G$.
\end{theorem}

For instance, in the odd wheel $K_1 \vee C_m$ with $m \ge 5$, the two triangles over a pair of consecutive cycle edges cover four vertices; for this pair Theorem~\ref{thm:selected} gives the upper bound $3$, whereas the global bound $n - k + 3$ equals $m$. For $X \subseteq V(G)$, let $G[X]$ denote the subgraph of $G$ induced by $X$. Note also that whenever $S_{\mathcal{A}}$ is a proper subset of $V(G)$, $G[S_{\mathcal{A}}]$ is $(k-1)$-colorable rather than $k$-critical, so the bound on $|\mathcal{A}|$ cannot be obtained by applying the global bound to $G[S_{\mathcal{A}}]$.

To prove Theorem~\ref{thm:selected}, we use colorings of edge deletions to locate vertices of degree $k-1$ in the selected cliques and then count incidences with these vertices.

\section{Preliminaries}\label{sec:prelim}

Stiebitz and Wessel~\cite{StiebitzWessel1993} characterized the critical graphs obtained by adding edges between a disjoint complete graph and an odd cycle. We use the following subgraph rigidity result of Stiebitz~\cite{Stiebitz1987}, in the form stated by Gao and Ma~\cite{GaoMa2023}, to recognize the extremal join as soon as it appears as an ordinary subgraph.

\begin{lemma}[Stiebitz~\cite{Stiebitz1987}]\label{lem:join}
Let $k \ge 4$ and $\ell \ge 3$. If a $k$-critical graph $G$ contains a subgraph isomorphic to $K_{k-3} \vee C_\ell$, then $\ell$ is odd and $G \cong K_{k-3} \vee C_\ell$.
\end{lemma}

Our second lemma combines clique-incidence vectors over $\QQ$ with $(k-1)$-colorings of edge-deleted subgraphs to locate vertices of degree $k - 1$ in the ambient graph $G$ when a family of cliques reaches the proposed bound. We write $N_G(v)$ for the set of neighbors of $v$ in $G$ and $d_G(v)=|N_G(v)|$.

\begin{lemma}\label{lem:degree}
Let $G$ be an $n$-vertex $k$-critical graph with $n > k \ge 3$, and let $\mathcal{A}$ be a nonempty family of distinct $q$-cliques of $G$, where $q = k - 1$. Let $N = |S_{\mathcal{A}}|$ and $t = |\mathcal{A}|$. If $t \ge N - q + 2$, then every clique in $\mathcal{A}$ contains at least two vertices of degree $q$ in $G$.
\end{lemma}

\begin{proof}
Write $S = S_{\mathcal{A}}$, and let $\QQ^S$ denote the space of rational vectors indexed by $S$. For each $R \in \mathcal{A}$, let $\mathbf{1}_R \in \QQ^S$ have coordinate $1$ on $R$ and $0$ on $S \setminus R$, and let $U$ be the subspace of $\QQ^S$ spanned by $\{\mathbf{1}_R : R \in \mathcal{A}\}$.

Let $(a_R)_{R \in \mathcal{A}}$ be rational numbers satisfying $\sum_{R \in \mathcal{A}} a_R = 0$, and write $z = \sum_{R \in \mathcal{A}} a_R \mathbf{1}_R$, with coordinate $z_v$ at $v \in S$. For each edge $xy \in E(G)$, write $\Lambda_a(xy) = \sum_{R \in \mathcal{A} : x,y \in R} a_R$, and fix a $q$-coloring of $G - xy$, the graph obtained by deleting the edge $xy$. Such a coloring exists because $G$ is $k$-critical. In this coloring, $x$ and $y$ receive the same color, since otherwise the coloring would be a $q$-coloring of $G$; let $C_{xy}$ be the color class containing $x$ and $y$. Every clique $R \in \mathcal{A}$ is contained in $S$. As observed by Gao and Ma~\cite{GaoMa2023}, if $x, y \in R$, then every vertex of $R \setminus \{x, y\}$ is adjacent to $x$ in $G - xy$, so $R \cap C_{xy} = \{x, y\}$ has size $2$; otherwise, $R$ remains a $q$-clique in $G - xy$, so its $q$ vertices receive distinct colors and $|R \cap C_{xy}| = 1$. Since $\sum_{R \in \mathcal{A}} a_R = 0$, summing $z_v$ over $C_{xy} \cap S$ gives
\begin{equation}\label{eq:common-color}
\sum_{v \in C_{xy} \cap S} z_v = \sum_{R \in \mathcal{A}} |R \cap C_{xy}|\, a_R = \sum_{R \in \mathcal{A}} a_R + \Lambda_a(xy) = \Lambda_a(xy).
\end{equation}

\textbf{Rational independence.}
We claim that $\dim U = t$. Suppose $\sum_{R \in \mathcal{A}} a_R \mathbf{1}_R = 0$ with $a_R \in \QQ$. Summing all coordinates over $S$ gives $q \sum_{R \in \mathcal{A}} a_R = 0$, so $\sum_{R \in \mathcal{A}} a_R = 0$ and \eqref{eq:common-color} gives $\Lambda_a(xy) = 0$ for every $xy \in E(G)$. For any $x \in S$, consider the restricted sum $\sum_{R \in \mathcal{A} : x \in R} a_R \mathbf{1}_R$. Its coordinate at $x$ equals the $x$-coordinate of $\sum_{R \in \mathcal{A}} a_R \mathbf{1}_R$, which is $0$; at every $y \in S \cap N_G(x)$ its coordinate is $\Lambda_a(xy) = 0$; and at every $y \in S \setminus (\{x\} \cup N_G(x))$ its coordinate is $0$ because no clique in $\mathcal{A}$ contains two nonadjacent vertices. Thus replacing $a_R$ by $0$ whenever $x \notin R$ yields another zero relation in $\QQ^S$. Fix any $Q \in \mathcal{A}$ and apply this restriction successively along the $q$ vertices of $Q$, each time starting from the current zero relation. Because all cliques in $\mathcal{A}$ are distinct sets of size $q$, the only term surviving all $q$ restrictions is $a_Q \mathbf{1}_Q$, so $a_Q = 0$. Hence $\dim U = t$.

\textbf{Low-degree vertices in each clique.}
Now fix an arbitrary clique $Q \in \mathcal{A}$, and let $E_Q \subseteq \QQ^S$ be the subspace of vectors whose coordinates vanish on $S \setminus Q$. Since $\dim \QQ^S = N$, $\dim U = t$, and $\dim E_Q = q$, we have $\dim(U \cap E_Q) \ge t + q - N \ge 2$. The linear functional on $U \cap E_Q$ sending a vector to the sum of its coordinates therefore has a nonzero vector $z$ in its kernel. Write $z = \sum_{R \in \mathcal{A}} a_R \mathbf{1}_R$ with $a_R \in \QQ$. Since the sum of the coordinates of $z$ equals $q \sum_{R \in \mathcal{A}} a_R$, we have $\sum_{R \in \mathcal{A}} a_R = 0$, so \eqref{eq:common-color} holds for every $xy \in E(G)$. Moreover, because $z \in E_Q$ vanishes on $S \setminus Q$, the left-hand side of \eqref{eq:common-color} reduces to $\sum_{v \in C_{xy} \cap Q} z_v$.

Fix $x \in Q$. For every $y \in Q \setminus \{x\}$, we have $x, y \in Q$ and hence $C_{xy} \cap Q = \{x, y\}$, so \eqref{eq:common-color} gives $\Lambda_a(xy) = z_x + z_y$. For every neighbor $y \in N_G(x) \setminus Q$ (including any neighbor outside $S$), the clique $Q$ omits $y$ and remains intact in $G - xy$; since $x \in C_{xy} \cap Q$, we have $C_{xy} \cap Q = \{x\}$, and \eqref{eq:common-color} gives $\Lambda_a(xy) = z_x$. Summing $\Lambda_a(xy)$ over all neighbors of $x$ in $G$ and using $\sum_{y \in Q} z_y = 0$ yields
\[
  \sum_{y \in N_G(x)} \Lambda_a(xy)
  = d_G(x) z_x + \sum_{y \in Q \setminus \{x\}} z_y
  = (d_G(x) - 1) z_x .
\]
On the other hand, each clique $R \in \mathcal{A}$ containing $x$ contains exactly $q - 1$ neighbors of $x$ in $G$, so the same sum equals
\[
  \sum_{R \in \mathcal{A} : x \in R} (q - 1) a_R = (q - 1) z_x .
\]
Comparing the two expressions gives $(d_G(x) - q) z_x = 0$, so $d_G(x) = q$ whenever $z_x \neq 0$. Since $z \in E_Q$ is nonzero and $\sum_{x \in Q} z_x = 0$, at least two vertices $x \in Q$ satisfy $z_x \neq 0$, and each has degree $q$ in $G$.
\end{proof}

\section{Proofs of the main theorems}\label{sec:proofs}

\begin{proof}[Proof of Theorem~\ref{thm:selected}]
Let $q = k - 1$, $S = S_{\mathcal{A}}$, $N = |S|$, and $t = |\mathcal{A}|$.\par\nopagebreak

\textbf{The join and $3$-critical cases.}
Suppose first that $k = 3$, or that $k \ge 4$ and $G$ contains a subgraph isomorphic to $K_{k-3} \vee C_\ell$ for some $\ell \ge 3$. If $k = 3$, then $\chi(G) = 3$ implies that $G$ contains an odd cycle; since every odd cycle has chromatic number $3$, the $3$-criticality of $G$ forces $G \cong C_n = K_0 \vee C_n$. If $k \ge 4$, then Lemma~\ref{lem:join} gives $G \cong K_{k-3} \vee C_\ell$ with $\ell$ odd. Thus in either case $G \cong K_{k-3} \vee C_m$, where $m = n - k + 3$ is odd, and $n > k$ gives $m \ge 5$.

Since $m \ge 5$, the cycle $C_m$ is triangle-free, so each $q$-clique of $K_{k-3} \vee C_m$ consists of $V(K_{k-3})$ together with the endpoints of a unique edge of $C_m$. The family $\mathcal{A}$ therefore corresponds bijectively to a nonempty edge set $F \subseteq E(C_m)$ of size $t$. If $F = E(C_m)$, then $\mathcal{A}$ consists of all $q$-cliques of $G$, $S = V(G)$, and $t = m = N - k + 3$. If $F$ is a nonempty proper subset of $E(C_m)$, then the subgraph of $C_m$ formed by $F$ and its incident vertices is a forest with $t$ edges and hence covers at least $t + 1$ vertices of $C_m$. Including the $k - 3$ vertices of $K_{k-3}$, which belong to every clique in $\mathcal{A}$, gives $N \ge (k - 3) + (t + 1) = t + k - 2$, so $t \le N - k + 2 < N - k + 3$.

\textbf{The join-free case for $k \ge 4$.}
Assume now that $k \ge 4$ and $G$ contains no subgraph isomorphic to $K_{k-3} \vee C_\ell$ for any $\ell \ge 3$. We prove that $t < N - k + 3$, that is, $t \le N - q + 1$. If $N = q$, then $t = 1 < 2 = N - q + 2$. Suppose for a contradiction that $N > q$ and $t \ge N - q + 2$; in particular, $t \ge 3$. Let $L = \{v \in S : d_G(v) = q\}$. By Lemma~\ref{lem:degree}, every clique in $\mathcal{A}$ contains at least two vertices of $L$.

We claim that each vertex $v \in L$ belongs to at most two cliques in $\mathcal{A}$. Take $Q \in \mathcal{A}$ with $v \in Q$. Because $|Q \setminus \{v\}| = q - 1$ and $d_G(v) = q$, the vertex $v$ has a unique neighbor $x \in N_G(v) \setminus Q$. Any $q$-clique through $v$ other than $Q$ is contained in $\{v\} \cup N_G(v) = Q \cup \{x\}$ and therefore has the form $\{v, x\} \cup W$ for some subset $W \subseteq (Q \setminus \{v\}) \cap N_G(x)$ of size $q - 2$. Since $G$ is $k$-critical and $n > k = q + 1$, the proper subgraph $G[Q \cup \{x\}]$ is $(k - 1)$-colorable and hence not complete. Thus $x$ is nonadjacent to at least one vertex of $Q \setminus \{v\}$, so $|(Q \setminus \{v\}) \cap N_G(x)| \le q - 2$. This leaves at most one choice for $W$, proving the claim.

Let $\ell = |L|$, and let $I$ be the number of pairs $(v, Q) \in L \times \mathcal{A}$ with $v \in Q$. By the preceding bounds, $2t \le I \le 2\ell$, so $\ell \ge t$. Each of the $N - \ell$ vertices in $D = S \setminus L$ belongs to at most $t$ cliques in $\mathcal{A}$, so counting all pairs $(v, Q) \in S \times \mathcal{A}$ with $v \in Q$ gives
\[
  qt \le I + (N - \ell)t \le 2\ell + (N - \ell)t .
\]
Rearranging gives $(t - 2)\ell \le (N - q)t$. Since $t \ge 3$ and $\ell \ge t$, we obtain $(t - 2)t \le (t - 2)\ell \le (N - q)t$, and hence $t \le N - q + 2$.

Since $t \ge N - q + 2$ by assumption, equality holds throughout: $t = N - q + 2$, $\ell = t$, $I = 2t = 2\ell$, and $qt = I + (N - \ell)t$. Consequently, every clique in $\mathcal{A}$ contains exactly two vertices of $L$, every vertex in $L$ belongs to exactly two cliques in $\mathcal{A}$, and every vertex in $D$ belongs to all $t$ cliques in $\mathcal{A}$. Because $\mathcal{A}$ is nonempty and each of its $q$-cliques consists of $D$ together with two vertices of $L$, the set $D$ is a clique of size $q - 2 = k - 3 \ge 1$ and every vertex of $L$ is adjacent to all of $D$. For each $v \in L$, the two distinct cliques in $\mathcal{A}$ containing $v$ have the form $D \cup \{v, u_1\}$ and $D \cup \{v, u_2\}$ with distinct $u_1, u_2 \in L \setminus \{v\}$. Together with $D$, the vertices $u_1$ and $u_2$ account for all $q$ neighbors of $v$ in $G$, so every vertex of $L$ has degree $2$ in $G[L]$. Since $\ell = t \ge 3$, the graph $G[L]$ contains a cycle $C_j$ for some $j \ge 3$. Together with $D$, this cycle forms a subgraph of $G$ isomorphic to $K_{k-3} \vee C_j$, contradicting our assumption.
\end{proof}

\begin{proof}[Proof of Theorem~\ref{thm:equality}]
Suppose first that $t_{k-1}(G) = m$. Since $n > k$, we have $m = n - k + 3 \ge 4$, so the family $\mathcal{A}$ of all $(k-1)$-cliques of $G$ is nonempty. Applying Theorem~\ref{thm:selected} to $\mathcal{A}$ gives $m = |\mathcal{A}| \le |S_{\mathcal{A}}| - k + 3 \le n - k + 3 = m$. Equality therefore holds in $|\mathcal{A}| \le |S_{\mathcal{A}}| - k + 3$, so Theorem~\ref{thm:selected} implies that $m$ is odd, $m \ge 5$, and $G \cong K_{k-3} \vee C_m$.

Conversely, suppose that $m \ge 5$ is odd and consider $H = K_r \vee C_m$, where $r = k - 3 \ge 1$. As observed in Section~\ref{sec:intro}, $t_{k-1}(H) = m$. It remains to verify that $H$ is $k$-critical. In any proper coloring of $H$, the vertices of $K_r$ receive $r$ distinct colors disjoint from those on $C_m$, and the odd cycle $C_m$ requires $3$ colors; combining an $r$-coloring of $K_r$ with a disjoint $3$-coloring of $C_m$ shows that $\chi(H) = r + 3 = k$. Now let $e \in E(H)$. If $e \in E(C_m)$, then $C_m - e$ is bipartite and $2$-colorable, so $H - e$ is $(r + 2)$-colorable. If $r \ge 2$ and $e \in E(K_r)$, assigning a common color to the endpoints of $e$, $r - 2$ new colors to the remaining vertices of $K_r$, and $3$ further colors to $C_m$ gives a $(k - 1)$-coloring of $H - e$. If $e = av$ with $a \in V(K_r)$ and $v \in V(C_m)$, assigning one color to both $a$ and $v$, $r - 1$ new colors to $V(K_r) \setminus \{a\}$, and $2$ further colors to the bipartite path $C_m - v$ obtained by deleting $v$ from $C_m$ again yields a $(k - 1)$-coloring of $H - e$. Thus $H - e$ is $(k - 1)$-colorable for every $e \in E(H)$. Since $H$ has no isolated vertices, every proper subgraph of $H$ omits at least one edge and is contained in $H - e$ for some $e \in E(H)$, so every proper subgraph of $H$ is $(k - 1)$-colorable. Hence $H$ is $k$-critical.
\end{proof}

\section*{Declaration of competing interest}

The author declares that he has no known competing financial interests or personal relationships that could have appeared to influence the work reported in this paper.

\section*{Acknowledgements}

The research was partially supported by RGC Competitive Earmarked Research Grant 16308821.

\section*{Data availability}

The proofs are contained in the article and do not rely on computational data or computer-generated certificates.

\section*{Declaration of generative AI use}

In the proof of Lemma~\ref{lem:degree}, ChatGPT assisted the author in developing the details of the common-color identity and the rational independence argument. It also assisted with language editing. The author independently checked all mathematical arguments and takes full responsibility for the accuracy and content of this work.

\end{document}